\documentclass{amsart}
\usepackage{amsmath,amsfonts,amsthm,amssymb,amscd}

\newtheorem{pro}{Proposition}[section]
\newtheorem{thm}{Theorem}[section]
\newtheorem{lem}{Lemma}[section]

\newtheorem{defi}{Definition}[section]
\newtheorem{remark}{Remark}[section]

\begin{document}
\title[K\"ahler-Ricci solitons with harmonic Bochner tensor]{On rigidity of
gradient K\"ahler-Ricci solitons with harmonic Bochner tensor}
\author{Qiang Chen}
\author{Meng Zhu}
\address{Department of Mathematics\\
Lehigh University\\
Bethlehem, PA 18015}
\email{qic208@lehigh.edu \& mez206@lehigh.edu}

\let\thefootnote\relax

\begin{abstract}
In this paper, we prove that complete gradient steady K\"ahler-Ricci
solitons with harmonic Bochner tensor are necessarily K\"ahler-Ricci flat,
i.e., Calabi-Yau, and that complete gradient shrinking (or expanding)
K\"ahler-Ricci solitons with harmonic Bochner tensor must be isometric to a
quotient of $N^k\times \mathbb{C}^{n-k}$, where $N$ is a K\"ahler-Einstein
manifold with positive (or negative) scalar curvature.
\end{abstract}

\maketitle

\section{Introduction}

A complete Riemannian manifold $(M, g_{ij})$ is called a Ricci soliton, if
there is a vector field $X$ and a constant $\lambda$ such that
\begin{equation*}
R_{ij}+ \frac{1}{2}(\nabla_i X_j + \nabla_j X_i) = \frac{\lambda}{2}g_{ij},
\end{equation*}
where $\lambda >0$, $\lambda=0$, or $\lambda <0$ corresponds to shrinking,
steady or expanding soliton respectively. Moreover, a Ricci soliton is
called a gradient Ricci soliton if the vector field is a gradient vector
field, i.e. $X=\nabla f$ for some smooth function $f$ on $M$. In this case,
the Ricci soliton equation above becomes
\begin{equation*}
R_{ij}+ \nabla_i\nabla_j f = \frac{\lambda}{2}g_{ij}.
\end{equation*}
Since R. Hamilton \cite{Ham1988} introduced the concept of Ricci solitons in
the mid 1980's, the study of Ricci solitons has attracted a lot of
attention. Ricci solitons, as self-similar solutions to Hamilton's Ricci
flow, are natural generalizations of Einstein metrics. Since they often arise
as the dilation limits of the singularities of the Ricci flow, Ricci
solitons play an important role in the singularity analysis of the Ricci
flow.\newline

In recent years, significant amount of results has been obtained in
classifying or understanding the geometry of shrinking solitons. In
particular, the classification of gradient shrinkers is known in dimensions
2 and 3, and assuming locally conformally flatness in all dimensions $n\ge 4$
, see e.g., \cite{P2}, \cite{NiWa2008}, \cite{CCZ2008}, \cite{ELM}, \cite
{PW3}, \cite{ Zha2009}, \cite{CWZ08}, \cite{Nab2010}, \cite{MuSe2009}, etc.
We also refer the reader to a recent survey paper of H.-D. Cao \cite{Cao2011} on
Ricci shrinkers. \newline

Regarding steady solitons, it is well-known that compact ones (as well as
compact expanding solitons) must be Einstein, see e.g., \cite{CaZh2006} for
a proof. However, much less is known for noncompact steady Ricci solitons.
In dimension $n=2$, R. Hamilton showed that any $2$-dimensional steady soliton
is isometric to the Cigar soliton, up to scaling. For $n=3$, G. Perelman \cite
{P2} conjectured that any 3-dimensional complete ($\kappa$-noncollapsed)
steady soliton with positive sectional curvature is isometric to the Bryant
soliton, which is the unique rotationally symmetric example on $\mathbb{R}^3$.
Very recently, H.-D. Cao and the first author \cite{CaCh2009} made the first
important progress on this problem. They showed that any $n$-dimensional, $
n\ge 3$, complete locally conformally flat steady Ricci soliton is either
flat or isometric to the Bryant soliton (for $n\ge 4$, Catino and Mantegazza
\cite{CaMa2010} independently proved this result by using a different
method). Subsequently, their work has been used by S. Brendle \cite{Bre2010} in
classifying 3-dimensional steady Ricci soliton satisfying certain asymptotic
condition, and X. X. Chen and Y. Wang \cite{ChWa2011} in classifying
4-dimensional half-comformally flat steady solitons, respectively. \newline

For more thorough discussions and results in Ricci solitons, the readers can
refer to the survey \cite{Cao2010} of H.-D. Cao and the references therein.
\newline

In this paper, we are interested in K\"{a}hler-Ricci solitons.

\begin{defi}
An n-dimensional K\"ahler manifold $(M^n, g_{i\bar{j}})$ is called a
gradient K\"ahler-Ricci soliton if there is a real-valued smooth function $f$
satisfying the soliton equation
\begin{equation*}
R_{i\bar{j}} + \nabla_i\nabla_{\bar{j}}f = \lambda g_{i\bar{j}},
\end{equation*}
for some constant $\lambda\in \mathbb{R}$ and such that $\nabla f$ is a
holomorphic vector field, i.e. $\nabla_i\nabla_j f=0$.
\end{defi}

In \cite{CaHa2000}, H.-D. Cao and R. Hamilton observed that complete noncompact
gradient steady K\"ahler-Ricci solitons with positive Ricci curvature such
that the scalar curvature attains the maximum must be Stein (and
diffeomorphic to $\mathbb{R}^{2n}$). Later, under the same assumptions, A. Chau
and L.-F. Tam \cite{ChTa2005}, and R. Bryant \cite{Bry2008} independently proved that
such steady K\"ahler-Ricci solitons are actually biholomorphic to $\mathbb{C}
^n$. Moreover, Chau and Tam showed that complete noncompact expanding K\"ahler-Ricci
solitons with nonnegative Ricci curvature are also biholomorphic to $\mathbb{
C}^n$.\newline

To state our results, let us first recall that on K\"ahler manifolds there
is a similar notion as the Weyl tensor, called the Bochner tensor. The
Bochner tensor $W_{i\bar{j}k\bar{l}}$ is defined by
\begin{eqnarray*}
W_{i\bar{j}k\bar{l}}&=&R_{i\bar{j}k\bar{l}} +\frac{R}{(n+1)(n+2)}(g_{i\bar{j}
}g_{k\bar{l}}+g_{i\bar{l}}g_{k\bar{j}}) \\
& & - \frac{1}{(n+2)}(R_{i\bar{j}}g_{k\bar{l}}+ R_{k\bar{l}}g_{i\bar{j}}+
R_{i\bar{l}}g_{k\bar{j}}+ R_{k\bar{j}}g_{i\bar{l}}),
\end{eqnarray*}
where $R_{i\bar{j}}=g^{k\bar{l}}R_{i\bar{j}k\bar{l}}$.
We also denote the divergence of the Bochner tensor by
\begin{eqnarray*}
C_{i\bar{j}k} &=& g^{q\bar{l}}\nabla_l W_{i\bar{j}k\bar{q}} \\
& = & \frac{n}{n+2}\nabla_i R_{k\bar{j}} - \frac{n}{(n+1)(n+2)}(g_{k\bar{j}}\nabla_i R
+ g_{i\bar{j}}\nabla_k R).
\end{eqnarray*}

\begin{defi}
\label{de2} A K\"ahler manifold $M^n$ is said to have harmonic Bochner
tensor if $C_{i\bar{j}k}=0$, i.e.,
\begin{equation*}
\nabla_i R_{k\bar{j}} = \frac{1}{n+1}(g_{k\bar{j}}\nabla_i R + g_{i\bar{j}}\nabla_k R
).
\end{equation*}
\end{defi}

%\begin{defi} \label{de2}
%Define
%\begin{eqnarray*}
%C_{i\bar{j}k} &=& \nabla_l W_{i\bar{j}k\bar{l}}\\
%& = & \frac{n}{n+2}\nabla_i R_{k\bar{j}} - \frac{n}{(n+1)(n+2)}(\nabla_i R g_{k\bar{j}} + \nabla_k R g_{i\bar{j}}).
%\end{eqnarray*}

%\hspace{-.4cm}If $C_{i\bar{j}k}=0$, we say that $M$ has harmonic Bochner tensor.
%\end{defi}

Very recently, by using a similar argument as in the paper \cite{CaCh2009} of H.-D. Cao and the first author, Y. Su
and K. Zhang \cite{SuZh2011} have shown that any complete noncompact
gradient K\"ahler-Ricci soliton with vanishing Bochner tensor is necessarily
K\"ahler-Einstein, and hence a quotient of the corresponding complex space
form. \newline

In this paper we investigate gradient K\"ahler-Ricci solitons with harmonic
Bochner tensor, and extend the classification result of Su and Zhang. Our main
results are:

\begin{thm}
\label{thm1} Any complete gradient steady K\"ahler-Ricci soliton with
harmonic Bochner tensor must be K\"ahler-Ricci flat (i.e., Calabi-Yau).
\end{thm}

\begin{thm}
\label{thm2} Any complete gradient shrinking (or expanding) K\"ahler-Ricci
soliton with harmonic Bochner tensor must be isometric to the quotient of $
N^k\times \mathbb{C}^{n-k}$, where $N^k$ is a $k$-dimensional
K\"ahler-Einstein manifold with positive (or negative) scalar curvature.
\end{thm}

\begin{remark}
It is known that a compact K\"ahler manifold with vanishing Bochner tensor (also called
Bochner-K\"ahler or Bochner-flat) is necessarily a compact quotient of $
M^{k}_{c}\times M^{n-k}_{-c}$, where $M^{k}_{c}$ and $M^{n-k}_{-c}$ denote
the complex space forms of constant holomorphic sectional curvature $c$ and $
-c$ respectively (cf., e.g., corollary 4.17 in \cite{Bry2001}). It follows
immediately that any compact K\"ahler-Ricci soliton with vanishing Bochner
tensor must be a quotient of complex space form.
\end{remark}

\begin{remark}
In the Riemanian case, by using a rigidity result of Petersen and Wylie \cite
{PeWy2009}, Fern\'{a}ndez-L\'{o}pez and Garc\'{\i}a-R\'{\i}o \cite{FeGa2010},
and Munteanu and Sesum \cite{MuSe2009} proved that Ricci shrinkers with
harmonic Weyl tensor must be rigid, i.e., a quotient of the product of an
Einstein manifold and $\mathbb{R}^k$.
\end{remark}

\hspace{-.4cm}\textbf{Acknowledgement} The authors would like to thank their
advisor, Professor Huai-Dong Cao, for suggesting this problem and for
stimulating discussions. We are also grateful for his constant encouragement
and support.

\section{proof of the main theorems}

Let $(M^n,g_{i\bar{j}},f)$ be a gradient K\"ahler-Ricci soliton, i.e.
\begin{equation*}
R_{i\bar{j}} + \nabla_i\nabla_{\bar{j}}f = \lambda g_{i\bar{j}}, \quad \text{
and}\quad \nabla_i\nabla_j f =0.
\end{equation*}

It is well-known that the following basic identities hold (see e.g. \cite
{CaZh2006}).

\begin{lem}
On a gradient K\"ahler-Ricci soliton, we have

\begin{equation}  \label{eq12}
R+|\nabla f|^2 - \lambda f = Const;
\end{equation}

\begin{equation}  \label{eq13}
R+ \Delta f = n\lambda;
\end{equation}

\begin{equation}  \label{eq14}
\nabla_i R_{k\bar{j}} = R_{i\bar{j}k\bar{l}}\nabla_l f ;
\end{equation}
and
\begin{equation}  \label{eq15}
\nabla_i R = R_{i\bar{j}}\nabla_j f.
\end{equation}
\end{lem}

From now on, we assume $(M^n,g_{i\bar{j}},f)$ is a gradient K\"ahler-Ricci
soliton with harmonic Bochner tensor so that
\begin{equation}  \label{eq1}
\nabla_i R_{k\bar{j}} = \frac{1}{n+1}(\nabla_i R g_{k\bar{j}} + \nabla_k R
g_{i\bar{j}}).
\end{equation}

\begin{lem}
\label{lem1} We have

\begin{eqnarray}  \label{eq2}
& & \lambda R_{i\bar{j}} - R_{i\bar{j}k\bar{l}}R_{\bar{k}l}  \notag \\
&=& \frac{1}{n+1}[ \frac{1}{n+1}\nabla_k R \nabla_{\bar{k}} f g_{i\bar{j}} +
(\lambda R-|Rc|^2)g_{i\bar{j}} - \frac{n}{n+1}\nabla_i R \nabla_{\bar{j}}f \\
& & \hspace{1cm} + \lambda R_{i\bar{j}} - R_{i\bar{k}}R_{k\bar{j}}],  \notag
\end{eqnarray}
and

\begin{eqnarray}  \label{eq3}
& & 2(n+1)\lambda\nabla_i R - 2R\nabla_i R - 2R_{i\bar{j}}\nabla_j R  \notag
\\
&=& -\frac{1}{n+1}\nabla_i R |\nabla f|^2 - \frac{1}{n+1}\nabla_k R\nabla_{
\bar{k}}f \nabla_i f.
\end{eqnarray}
\end{lem}

\medskip

\noindent \textbf{Proof:} On one hand, by differentiating (2.4), we obtain

\begin{equation*}
\Delta R = \nabla_{k}\nabla_{\bar{k}}R = \nabla_k R \nabla_{\bar{k}}f + R_{k
\bar{l}}\nabla_{\bar{k}}\nabla_l f.
\end{equation*}

From (\ref{eq1}), we get

\begin{eqnarray}  \label{eq4}
\nabla_k\nabla_{\bar{k}}R_{i\bar{j}} &=& \frac{1}{n+1}(\Delta R g_{i\bar{j}}
+ \nabla_i\nabla_{\bar{j}} R)  \notag \\
&=& \frac{1}{n+1}(\nabla_k R \nabla_{\bar{k}} f g_{i\bar{j}} + R_{k\bar{l}
}\nabla_{\bar{k}}\nabla_l f g_{i\bar{j}} + \nabla_i R_{k\bar{j}}\nabla_{\bar{
k}} f + R_{k\bar{j}}\nabla_i\nabla_{\bar{k}}f)  \notag \\
&=& \frac{1}{n+1}[\nabla_k R \nabla_{\bar{k}} f g_{i\bar{j}} + (\lambda
R-|Rc|^2)g_{i\bar{j}} + \frac{1}{n+1}\nabla_i R \nabla_{\bar{j}}f  \notag \\
& & \hspace{1cm} + \frac{1}{n+1} \nabla_k R \nabla_{\bar{k}} f g_{i\bar{j}}
+ \lambda R_{i\bar{j}} - R_{i\bar{k}}R_{k\bar{j}}]. \\
\notag
\end{eqnarray}

On the other hand, by differentiating (2.3), we have
% since $$\nabla_{\bar{k}}R_{i\bar{j}} = R_{\bar{k}i\bar{j}l}\nabla_{\bar{l}}f,$$

\begin{eqnarray*}
\nabla_k\nabla_{\bar{k}} R_{i\bar{j}} &=& \nabla_iR_{\bar{j}l}\nabla_{\bar{l}
}f + R_{i\bar{j}k\bar{l}}\nabla_{\bar{k}}\nabla_{l}f \\
&=& \nabla_k R_{i\bar{j}}\nabla_{\bar{k}}f + R_{i\bar{j}k\bar{l}}\nabla_{
\bar{k}}\nabla_l f \\
&=& \nabla_k R_{i\bar{j}}\nabla_{\bar{k}}f + \lambda R_{i\bar{j}} - R_{i\bar{
j}k\bar{l}}R_{\bar{k}l}.
\end{eqnarray*}

Now, by plugging in formula (\ref{eq4}), we obtain (\ref{eq2}).\newline

Next, by taking the divergence on both sides of (\ref{eq2}), we get

\begin{eqnarray*}
& & \lambda\nabla_i R - (\nabla_i R_{k\bar{l}}) R_{\bar{k}l} - R_{i\bar{j}k
\bar{l}}\nabla_{j}R_{\bar{k}l} \\
&=& \frac{1}{n+1}[ \frac{1}{n+1}\nabla_i\nabla_k R \nabla_{\bar{k}} f +
\frac{1}{n+1}\nabla_k R \nabla_i\nabla_{\bar{k}} f + \lambda\nabla_i R
-\nabla_i|Rc|^2 \\
& & \hspace{1cm} - \frac{n}{n+1}\nabla_j\nabla_i R \nabla_{\bar{j}}f - \frac{
n}{n+1}\nabla_i R\Delta f + \lambda \nabla_i R - (\nabla_j R_{i\bar{k}})R_{k
\bar{j}} - R_{i\bar{k}}\nabla_k R] \\
&=& \frac{1}{n+1}[ \frac{1}{n+1}\nabla_iR_{k\bar{l}}\nabla_{l}f \nabla_{\bar{
k}} f + \frac{\lambda}{n+1}\nabla_i R -\frac{1}{n+1}R_{i\bar{k}}\nabla_k R +
\lambda\nabla_i R -2 R_{k\bar{l}}\nabla_i R_{\bar{k}l} \\
& & \hspace{1cm} - \frac{n}{n+1}\nabla_i R_{j\bar{k}}\nabla_k f \nabla_{\bar{
j}}f - \frac{\lambda n^2}{n+1}\nabla_i R + \frac{n}{n+1}R\nabla_i R +
\lambda \nabla_i R \\
& & \hspace{1cm} - R_{k\bar{j}}\nabla_i R_{j\bar{k}} - R_{i\bar{k}}\nabla_k
R] \\
\end{eqnarray*}

That is,

\begin{eqnarray*}
& & \lambda\nabla_i R - (\nabla_i R_{k\bar{l}}) R_{\bar{k}l} - R_{i\bar{j}k
\bar{l}}\nabla_{j}R_{\bar{k}l}  \notag \\
&=& \frac{1}{n+1}[ -\frac{n-1}{(n+1)^2}\nabla_iR|\nabla f|^2 - \frac{n-1}{
(n+1)^2}\nabla_k R\nabla_{\bar{k}}f \nabla_i f \\
& & \hspace{1cm} + (3-n)\lambda\nabla_i R -(1+\frac{1}{n+1})R_{i\bar{k}
}\nabla_k R -3 R_{k\bar{l}}\nabla_i R_{\bar{k}l}+ \frac{n}{n+1}R\nabla_i R].
\notag
\end{eqnarray*}

But,

\begin{eqnarray*}
R_{l\bar{k}}\nabla_i R_{k\bar{l}} &=& \frac{1}{n+1}R_{l\bar{k}}(\nabla_i
Rg_{k\bar{l}}+ \nabla_k Rg_{i\bar{l}}) \\
&=& \frac{1}{n+1} R\nabla_i R + \frac{1}{n+1}R_{i\bar{j}}\nabla_j R,
\end{eqnarray*}

\hspace{-.4cm}and

\begin{eqnarray*}
R_{i\bar{j}k\bar{l}}\nabla_j R_{l\bar{k}} &=& \frac{1}{n+1}R_{i\bar{j}k\bar{l
}}(\nabla_j Rg_{l\bar{k}}+ \nabla_l Rg_{j\bar{k}}) \\
&=& \frac{1}{n+1}R_{i\bar{j}}\nabla_j R + \frac{1}{n+1} R_{i\bar{l}}\nabla_l
R \\
&=& \frac{2}{n+1}R_{i\bar{j}}\nabla_j R.
\end{eqnarray*}

\hspace{-.4cm}Hence, we have

\begin{eqnarray*}
& &\lambda\nabla_i R - \frac{1}{n+1} R\nabla_i R - \frac{3}{n+1}R_{i\bar{j}
}\nabla_j R \\
&=&\lambda\nabla_i R - (\nabla_i R_{k\bar{l}}) R_{\bar{k}l} - R_{i\bar{j}k
\bar{l}}\nabla_{j}R_{\bar{k}l} \\
&=& \frac{1}{n+1}[ -\frac{n-1}{(n+1)^2}\nabla_iR|\nabla f|^2 - \frac{n-1}{
(n+1)^2}\nabla_k R\nabla_{\bar{k}}f \nabla_i f \\
& & \hspace{1cm} + (3-n)\lambda\nabla_i R -(1+\frac{1}{n+1})R_{i\bar{k}
}\nabla_k R -3 R_{k\bar{l}}\nabla_i R_{\bar{k}l}+ \frac{n}{n+1}R\nabla_i R].
\end{eqnarray*}

Therefore, formula (\ref{eq3}) follows easily. \qed\newline

Now, suppose that $\nabla f \neq 0$ at some point $p$. Then we may choose an
orthonormal frame $\{e_1,e_2, \cdots, e_n\}$ of holomorphic vector fields at
$p$ such that $e_1$ is parallel to $\nabla f$. Therefore, we have $|\nabla_1
f|= |\nabla f|$ and $\nabla_k f = 0$ for $k=2,\cdots, n$.

\begin{lem}
\label{lem2}

Suppose $\nabla f\neq 0$ at $p$. Then, under the frame $\{e_1,e_2, \cdots,
e_n\}$ chosen above, we have

%\begin{equation} \label{eq5}
%R_{i\bar{j}k\bar{1}} = \frac{1}{n+1}(R_{i\bar{1}}g_{k\bar{j}} + R_{k\bar{1}}g_{i\bar{j}}).
%\end{equation}

\begin{equation*}
R_{k\bar{1}}=R_{1\bar{k}}=0 \quad \mathit{for} \quad k\geq 2.
\end{equation*}
\end{lem}

\hspace{-.4cm}\textbf{Proof:} From (\ref{eq14}) and (\ref{eq1}), we have at $
p$,
\begin{equation*}
R_{i\bar{j}k\bar{1}}\nabla_1 f = \frac{1}{n+1}(\nabla_i R g_{k\bar{j}} +
\nabla_k R g_{i\bar{j}}) = \frac{1}{n+1}(R_{i\bar{1}}g_{k\bar{j}} + R_{k\bar{
1}}g_{i\bar{j}})\nabla_1 f.
\end{equation*}

\hspace{-.4cm}It follows that
\begin{equation*}
R_{i\bar{j}k\bar{1}} = \frac{1}{n+1}(R_{i\bar{1}}g_{k\bar{j}} + R_{k\bar{1}
}g_{i\bar{j}}).
\end{equation*}
In particular, for $k\geq 2$, we have that
\begin{equation*}
R_{1\bar{1}k\bar{1}} = \frac{1}{n+1}R_{k\bar{1}} \quad \textrm{and} \quad R_{1\bar{k}1\bar{1}}=0.
\end{equation*}
However, on the other hand, it is easy to see that
\begin{equation*}
R_{1\bar{1}k\bar{1}} = \overline{R_{\bar{1}1\bar{k}1}} = \overline{R_{1\bar{k
}1\bar{1}}}= 0.
\end{equation*}
Therefore, $R_{k\bar{1}}=R_{1\bar{k}}=0$ for $k\geq 2$. \qed \newline

Lemma \ref{lem2} tells us that $\nabla f$ is an eigenvector of the Ricci
curvature tensor. Thus we may choose another orthonormal frame $
\{w_{1}=e_{1},w_{2},\cdots ,w_{n}\}$ at $p$ such that  $|\nabla
_{1}f|=|\nabla f|$ \ and the Ricci curvature is diagonalized at $p$, i.e.
\begin{equation*}
R_{i\bar{j}}=R_{i\bar{i}}\delta _{ij}.
\end{equation*}

\begin{pro}
Suppose that $\nabla f \neq 0$ at $p$. Then under the orthonormal frame $
\{w_1,w_2,\cdots, w_n\}$ chosen above, we have the following identities at $p
$

\begin{equation}  \label{eq6}
n \lambda R_{1\bar{1}} - RR_{1\bar{1}} = \lambda R-|Rc|^2 - \frac{n-1}{n+1}
R_{1\bar{1}} |\nabla f|^2,
\end{equation}

\hspace{-.4cm}and

\begin{equation}  \label{eq7}
(n+1)\lambda R_{1\bar{1}} - RR_{1\bar{1}} - R_{1\bar{1}}^2 = -\frac{1}{n+1}
R_{1\bar{1}}|\nabla f|^2.
\end{equation}
\end{pro}

\hspace{-.4cm}\textbf{Proof:} In (\ref{eq2}), setting $i=j=1$, we have

\begin{eqnarray*}
&& \lambda R_{1\bar{1}} - \frac{1}{n+1}R_{1\bar{1}}^2 - \frac{1}{n+1}RR_{1
\bar{1}} \\
&=& \lambda R_{1\bar{1}} - \frac{2}{n+1}R_{1\bar{1}}^2 - \frac{1}{n+1}R_{1
\bar{1}}(R-R_{1\bar{1}}) \\
&=& \lambda R_{1\bar{1}} - \frac{2}{n+1}R^2_{1\bar{1}}- \frac{1}{n+1}R_{1
\bar{1}}\sum_{k=2}^n R_{k\bar{k}} \\
&=& \lambda R_{1\bar{1}} - R_{1\bar{1}1\bar{1}}R_{1\bar{1}} -
\sum_{k=2}^{n}R_{1\bar{1}k\bar{k}}R_{k\bar{k}} \\
&=& \lambda R_{1\bar{1}} - \sum_{k=1}^{n} R_{1\bar{1}k\bar{k}}R_{k\bar{k}} \\
&=& \frac{1}{n+1}[ \frac{1}{n+1}R_{1\bar{1}}|\nabla f|^2 + \lambda R-|Rc|^2
- \frac{n}{n+1}R_{1\bar{1}} |\nabla f|^2 + \lambda R_{1\bar{1}} - R_{1\bar{1}
}^2]. \\
\end{eqnarray*}

Thus, formula (\ref{eq6}) follows immediately.\newline

Next, by setting $i=1$ in (\ref{eq3}) and dividing both sides of the equation
by $\nabla_1 f$, we get (\ref{eq7}). \qed\newline

\begin{pro}
At a point $p$ where $\nabla f \neq 0$, we have either
\begin{equation*}
Rc(\nabla f, \nabla f)=0,
\end{equation*}
or
\begin{equation*}
Rc(\nabla f, \nabla f) = \frac{\lambda}{n+4} |\nabla f|^2.
\end{equation*}
\end{pro}

\hspace{-.4cm}\textbf{Proof:} Since at point $p$, $\nabla f\neq 0$, formula (
\ref{eq7}) implies that in a neighborhood of $p$ we have

\begin{equation}
\left[(n+1)\lambda - R - \frac{R_{j\bar{i}}\nabla_i f\nabla_{\bar{j}}f}{
|\nabla f|^2} + \frac{1}{n+1}|\nabla f|^2\right] \frac{R_{j\bar{i}}\nabla_i
f\nabla_{\bar{j}}f}{|\nabla f|^2} = 0.
\end{equation}

Therefore, there are two possibilities

\textbf{I)}\ $R_{j\bar{i}}\nabla_i f\nabla_{\bar{j}}f = 0$ at $p$,\newline

or\newline

\textbf{II)}\ $R_{j\bar{i}}\nabla_i f\nabla_{\bar{j}}f\neq 0$ at $p$. In
this case, near $p$ we have

\begin{equation*}
-(n+1)\lambda + R + \frac{R_{j\bar{i}}\nabla_i f\nabla_{\bar{j}}f}{|\nabla
f|^2} - \frac{1}{n+1}|\nabla f|^2=0.
\end{equation*}

Taking covariant derivative on both sides gives us

\begin{eqnarray*}
0 &=& \nabla_k R + \frac{1}{|\nabla f|^2}(\nabla_i f\nabla_{\bar{j}}
f\nabla_k R_{j\bar{i}} + R_{j\bar{i}}\nabla_i f\nabla_k\nabla_{\bar{j}} f) -
\frac{\nabla_j f\nabla_k\nabla_{\bar{j}} f}{|\nabla f|^4}R_{l\bar{i}
}\nabla_i f\nabla_{\bar{l}}f \\
& & - \frac{1}{n+1}(\nabla_j f\nabla_k\nabla_{\bar{j}}f) \\
&=& \nabla_k R + \frac{1}{(n+1)|\nabla f|^2}\nabla_i f\nabla_{\bar{j}}
f(\nabla_k R g_{j\bar{i}} + \nabla_j R g_{k\bar{i}}) + \frac{1}{|\nabla f|^2}
(\lambda\nabla_k R - R_{k\bar{j}}\nabla_j R) \\
& & - \frac{\lambda\nabla_k f - \nabla_k R}{|\nabla f|^4}\nabla_i R \nabla_{
\bar{i}}f - \frac{1}{n+1}(\lambda\nabla_k f - \nabla_k R).
\end{eqnarray*}

Evaluating the identity above at $p$ under the orthonormal frame $\{w_1,
w_2, \cdots, w_n\}$ yields

\begin{eqnarray*}
0 &=& R_{1\bar{1}} + \frac{2}{(n+1)|\nabla f|^2}R_{1\bar{1}}|\nabla f|^2 +
\frac{1}{|\nabla f|^2}(\lambda R_{1\bar{1}} - R^2_{1\bar{1}}) \\
& & - \frac{\lambda - R_{1\bar{1}}}{|\nabla f|^4}R_{1\bar{1}} |\nabla f|^2 -
\frac{1}{n+1}(\lambda - R_{1\bar{1}}) \\
& =& \frac{n+4}{n+1}R_{1\bar{1}} - \frac{1}{n+1}\lambda.
\end{eqnarray*}

Thus, we have $Rc(\nabla f, \nabla f)=\frac{\lambda}{n+4}|\nabla f|^2$
whenever $Rc(\nabla f, \nabla f)\neq 0$. \qed \newline

Now we are ready to prove the main theorems. \newline

First, we may assume that $f\neq Const$, for otherwise we get that $M$ is
K\"ahler-Einstein from the soliton equation.\newline

\hspace{-.4cm}\textbf{Proof of theorem \ref{thm1}:} For steady K\"{a}hler-Ricci solitons, we have $\lambda =0$.\newline

From proposition 2.2, we know that $Rc(\nabla f, \nabla f)=\frac{\lambda}{n+4
}|\nabla f|^2=0$ whenever $Rc(\nabla f, \nabla f)\neq 0$, which is a
contradiction. Therefore, we always have $Rc(\nabla f, \nabla f) = 0$. Then (\ref{eq6}) implies that $Rc = 0$ in the set $\{p\in M| \nabla f(p)\neq 0\}$.
On the other hand, by the soliton equation, it is easy to see that we also
have $Rc=0$ in the interior of the set $\{p\in M| \nabla f(p)=0\}$.
Thus the steady soliton $M$ must be K\"ahler-Ricci flat. \qed\newline

\hspace{-.4cm}\textbf{Proof of theorem \ref{thm2}:} For shrinking and
expanding K\"{a}hler-Ricci solitons, we have $\lambda \neq 0$.\newline

In this case, from proposition 2.2 and the continuity of $\frac{Rc(\nabla f,
\nabla f)}{|\nabla f|^2}$, we conclude that in each component of the open
set $A=\{p\in M| \nabla f(p)\neq 0\}$, we have either $Rc(\nabla f, \nabla f)=
\frac{\lambda}{n+4}|\nabla f|^2$ or $Rc(\nabla f, \nabla f) = 0$.\newline

If $Rc(\nabla f, \nabla f)=\frac{\lambda}{n+4}|\nabla f|^2$ in some component $\Omega$ of $A$, then at
any point $p\in \Omega$ we have $R_{1\bar{1}}=\frac{\lambda}{
n+4}$ and $\nabla R(p) = \frac{\lambda}{n+4}\nabla f(p)$. Therefore, we have
$\nabla R = \frac{\lambda}{n+4}\nabla f $ in $\Omega$. It then follows that $R
= \frac{\lambda}{n+4} f + C$ in $\Omega$. Thus (\ref{eq7}) implies that $
|\nabla f|^2 = \frac{n+1}{n+4}\lambda f + C$ in $\Omega$. Since $R+|\nabla f|^2 -
\lambda f = C_0$, we have $f = C_1$ in $\Omega$, which contradicts to the
fact that $\nabla f \neq 0$ in $\Omega$.\newline

Therefore, we must have $Rc(\nabla f, \nabla f)= 0$ in $A$. Since $f=Const$ in the interior of $
M\backslash A$, we have $Rc(\nabla f, \nabla f)= 0$ on the whole manifold $M$. It follows that $\nabla R =0$ on $M$. Then (\ref{eq1}) implies that the
Ricci curvature tensor is parallel on $M$. Therefore, by the de-Rahm
decomposition theorem, the universal cover of $M$ is isometric to $
N^{n-1}\times \mathbb{C}$, where $N$ is again an $n-1$ dimensional
K\"ahler-Ricci soliton with harmonic Bochner tensor. Thus by induction, we
can finally get that $M$ is isometric to a quotient of the product of a
K\"ahler-Einstein manifold and the complex Euclidean space.\newline

This finishes the proof. \qed

\end{document}